\documentclass[12pt,a4paper]{article}
\usepackage{pifont}
\usepackage{amsfonts}
\usepackage{mathrsfs}
\usepackage{amssymb}
\usepackage{amsmath}

\newtheorem{theorem}{Theorem}
\makeatletter

\newcommand{\Rmnum}[1]{\expandafter\@slowromancap\romannumeral #1@}
\makeatother

\begin{document}
\title{A simple proof of Jordan normal form}
\author{
 Yuqun
Chen  \\
{\small \ School of Mathematical Sciences, South China Normal
University}\\
{\small Guangzhou 510631, P.R. China}\\
{\small yqchen@scnu.edu.cn}}

\date{}

\maketitle \noindent\textbf{Abstract:} In this note, a simple proof
Jordan normal form and rational form of matrices over a field is
given.

\ \

Let $F$ be a field, $M_n(F)$ the set of $n\times n$ matrices over
$F$, $\lambda_i\in F,\ n_i$ a natural number, $i=1,2,\cdots,r$. A
 matrix of the form

\begin{eqnarray}\label{e1}
\left( \begin{array}{ccccc}
               J_{\lambda_1,n_1} & 0 & \dots & 0 \\
               0 & J_{\lambda_2,n_2} &  \dots & 0 \\
               \vdots & \vdots &  \vdots & \vdots \\
               0 & 0 &  \dots & J_{\lambda_r,n_r}
               \end{array} \right)
\end{eqnarray}
is called a Jordan normal form, where for each $i$,
\begin{eqnarray*}
J_{\lambda_i,n_i}=\left( \begin{array}{cccccc}
               \lambda_1 & 1&0 & \dots & 0 & 0\\
               0 & \lambda_1 &1&  \dots & 0 & 0 \\
               \vdots & \vdots &  \vdots & \vdots & \vdots & \vdots\\
               0 & 0 & 0& \dots & \lambda_1& 1 \\
               0 & 0 & 0 &\dots& 0 & \lambda_1
               \end{array} \right)
\end{eqnarray*}
is an $n_i\times n_i$ matrix.

Let $A,B\in M_n(F)$. Then $A$ is similar to $B$ in $M_n(F)$ if there
exists an invertible matrix $P\in M_n(F)$ such that $A=PBP^{-1}$.

It is known that an $n\times n$ matrix over the complex field is
similar to a Jordan normal form, see \cite{Ayres}. In this paper, a
short proof of this result is given.

\begin{theorem}\label{t1} (Jordan)
Let $F[\lambda]$ be the polynomial ring with one variable $\lambda$,
$A\in M_n(F)$ and $\lambda E-A$ the eigenmatrix. Suppose that there
exist invertible $\lambda$-matrices $P(\lambda),\ Q(\lambda)$ such
that
\begin{eqnarray}\label{e2}
P(\lambda)(\lambda E-A)Q(\lambda)=\left( \begin{array}{ccccccc}
               1 & \dots&0 & 0& \dots & 0 \\
               \vdots & \vdots &  \vdots & \vdots  & \vdots & \vdots  \\
               0 & \dots &1& 0&  \dots & 0 \\
               0 & \dots & 0& (\lambda-\lambda_1)^{n_1}   & \cdots & 0 \\
               \vdots & \vdots &  \vdots & \vdots  & \vdots & \vdots \\
               0 & \dots & 0& 0 &\dots & (\lambda-\lambda_r)^{n_r}
               \end{array} \right)
\end{eqnarray}
Then in $M_n(F)$, $A$ is similar to (\ref{e1}).
\end{theorem}
\noindent{\bf Proof:} Let $V$ be a $n$-dimensional vector space over
$F$ with a basis $e_1,\dots,e_n$. Let $\cal A$ be a linear
transformation of $V$ such that
\begin{eqnarray*}
({\cal A}E-A)\left( \begin{array}{cccccc}
               e_1 \\
                 \vdots  \\
               e_n
               \end{array} \right)=0.
\end{eqnarray*}
Since $Q(\lambda)$ is invertible, each entry in $Q({\cal A})^{-1}$
is a ${\cal A}$-polynomial. Now, by (\ref{e2}), we have
\begin{eqnarray*}
\left( \begin{array}{ccccccc}
               1 & \dots&0 & 0& \dots & 0 \\
               \vdots & \vdots &  \vdots & \vdots  & \vdots  & \vdots \\
               0 & \dots &1& 0&  \dots & 0 \\
               0 & \dots & 0& ({\cal A}-\lambda_1)^{n_1}   & \cdots & 0 \\
               \vdots & \vdots &  \vdots & \vdots  & \vdots & \vdots \\
               0 & \dots & 0& 0 &\dots & ({\cal A}-\lambda_r)^{n_r}
               \end{array} \right)
\left( \begin{array}{cc}
               \ast \\
                 \vdots  \\
               \ast\\
               y_1 \\
                 \vdots  \\
               y_r
               \end{array} \right)=0,
\end{eqnarray*}
where
\begin{eqnarray*}
\left( \begin{array}{cc}
               \ast \\
                 \vdots  \\
               \ast\\
               y_1 \\
                 \vdots  \\
               y_r
               \end{array} \right)
=Q({\cal A})^{-1} \left( \begin{array}{cccccc}
               e_1 \\
                 \vdots  \\
               e_n
               \end{array} \right).
\end{eqnarray*}
Then we have $({\cal A}-\lambda_1)^{n_1}y_1=0,\dots,({\cal
A}-\lambda_r)^{n_r}y_r=0$.

For any $\alpha\in V$, suppose that
\begin{eqnarray*}
\alpha=(a_1,\dots,a_n)\left( \begin{array}{cccccc}
               e_1 \\
                 \vdots  \\
               e_n
               \end{array} \right), \ a_i\in F,\ i=1,\dots,n.
\end{eqnarray*}
Then
\begin{eqnarray*}
\alpha &=&(a_1,\dots,a_n)Q({\cal A})\left( \begin{array}{cc}
               \ast \\
                 \vdots  \\
               \ast\\
               y_1 \\
                 \vdots  \\
               y_r
               \end{array} \right)
=(\star,\dots,\star,f_1({\cal A}),\dots,f_r({\cal A})) \left(
\begin{array}{cc}
               \ast \\
                 \vdots  \\
               \ast\\
               y_1 \\
                 \vdots  \\
               y_r
               \end{array} \right)\\
&=&f_1({\cal A})y_1+\dots+f_r({\cal A})y_r.
\end{eqnarray*}
Noting that $n_1+\dots+n_r=n$, it is easy to see that
$$
y_1,\ ({\cal A}-\lambda_1)y_1,\dots, ({\cal
A}-\lambda_1)^{n_1-1}y_1,\dots\dots,y_r,\ ({\cal
A}-\lambda_r)y_r,\dots, ({\cal A}-\lambda_r)^{n_r-1}y_r
$$
forms a basis of $V$ and the matrix of ${\cal A}$ under this basis
is (\ref{e1}).

The proof is completed.

\ \

The following known theorem can be similarly proved.

\begin{theorem} Let the notion be as in Theorem \ref{t1}.
Suppose that there exist invertible $\lambda$-matrices $P(\lambda),\
Q(\lambda)$ such that
\begin{eqnarray}\label{e4}
P(\lambda)(\lambda E-A)Q(\lambda)=\left( \begin{array}{ccccccc}
               1 & \dots&0 & 0& \dots & 0 \\
               \vdots & \vdots &  \vdots & \vdots  & \vdots  & \vdots \\
               0 & \dots &1& 0&  \dots & 0 \\
               0 & \dots & 0& g_1(\lambda)   & \cdots & 0 \\
               \vdots & \vdots &  \vdots & \vdots  & \vdots & \vdots \\
               0 & \dots & 0& 0 &\dots & g_r(\lambda)
               \end{array} \right)
\end{eqnarray}
where
$g_i(\lambda)=\lambda^{n_i}-a_{in_i-1}\lambda^{n_i-1}-\dots-a_{i1}\lambda-a_{i0}\in
F[\lambda],\ i=1,\dots,r$ and $n_1+\dots+n_r=n$.

Then in $M_n(F)$, $A$ is similar to a rational form
\begin{eqnarray}\label{e3}
\left( \begin{array}{ccccc}
               T_1 & 0 & \dots & 0 \\
               0 & T_2 &  \dots & 0 \\
               \vdots & \vdots &  \vdots & \vdots \\
               0 & 0 &  \dots & T_r
               \end{array} \right)
\end{eqnarray}
where for each $i$,
\begin{eqnarray*}
T_i=\left( \begin{array}{cccccc}
               0 & 1&0 & \dots & 0 & 0\\
               0 & 0 &1&  \dots & 0 & 0 \\
               \vdots & \vdots &  \vdots & \vdots & \vdots & \vdots \\
               0 & 0 & 0& \dots & 0& 1 \\
               a_{i0} & a_{i1} & a_{i2} &\dots& a_{in_i-2} & a_{in_i-1}
               \end{array} \right).
\end{eqnarray*}
\end{theorem}
\noindent{\bf Proof:} By the proof as the same as Theorem \ref{t1},
$$
y_1,\ {\cal A}y_1,\dots, {\cal A}^{n_1-1}y_1,\dots\dots,y_r,\ {\cal
A}y_r,\dots, {\cal A}^{n_r-1}y_r
$$
forms a basis of $V$ and the matrix of ${\cal A}$ under this basis
is (\ref{e3}). The proof is completed.

\ \

\noindent{\bf Remark:} If $F$ is  algebraically closed, then one can
find invertible $\lambda$-matrices $P(\lambda),\ Q(\lambda)$ such
that (\ref{e2}) holds. Generally, we have (\ref{e4}).

\end{document}